\documentclass{article}

\usepackage{arxiv}
\usepackage{authblk}
\usepackage{graphicx}      
\usepackage{natbib}        
\usepackage{amsmath}
\usepackage{amsfonts}
\usepackage{amssymb}
\usepackage{float}
\usepackage[ruled,vlined, linesnumbered, algo2e]{algorithm2e}
\usepackage{placeins}

\let\Oldsection\section
\renewcommand{\section}{\FloatBarrier\Oldsection}

\let\Oldsubsection\subsection
\renewcommand{\subsection}{\FloatBarrier\Oldsubsection}

\let\Oldsubsubsection\subsubsection
\renewcommand{\subsubsection}{\FloatBarrier\Oldsubsubsection}

\makeatletter
\newcommand{\removelatexerror}{\let\@latex@error\@gobble}
\makeatother

\newcommand{\state}{\mathbf{x}}
\newcommand{\control}{\mathbf{u}}
\newcommand{\uspace}{\mathcal{U}}
\newcommand{\target}{\mathcal{T}}
\newtheorem{remark}{Remark}
\begin{document}

\title{An Actor Critic Method for Free Terminal Time Optimal Control} 


\author[1]{Evan Burton} 
\author[1,2]{Tenavi Nakamura-Zimmerer} 
\author[1]{Qi Gong}
\author[3,1]{Wei Kang}

\affil[1]{University of California Santa Cruz, 
   Santa Cruz, CA 95060, USA (e-mail: esburton@ucsc.edu, tenakamu@ucsc.edu, qgong@ucsc.edu).}
\affil[2]{Flight Dynamics Branch, NASA Langley Research Center, Hampton, VA 23666, USA.
}
\affil[3]{Naval Postgraduate School, Monterey, CA 93943, USA
 (email: wkang@nps.edu)}

\maketitle

\begin{abstract}                
Optimal control problems with free terminal time present many challenges including nonsmooth and discontinuous control laws, irregular value functions, many local optima, and the curse of dimensionality. To overcome these issues, we propose an adaptation of the model-based actor-critic paradigm from the field of Reinforcement Learning via an exponential transformation to learn an approximate feedback control and value function pair. We demonstrate the algorithm's effectiveness on prototypical examples featuring each of the main pathological issues present in problems of this type.
\end{abstract}


\section{Introduction}

We aim to solve optimal feedback control problems with a free terminal time and fixed target state, which common in trajectory optimization path-planning. These problems often have solutions which saturate the control bounds and have value functions which may not even be continuous.  We will consider the problem of minimizing a running cost $\ell(\state, \control) \geq 1$ with free terminal time and given terminal state $x(t_f) \in \mathcal{T}$. That is, we wish to approximate the solution of
\begin{equation}
\begin{array}{c}
V\left( \state_0 \right) = \underset{\control \in \uspace, t_f \geq 0}{\min} \int_{0}^{t_f} \ell\left(\state\left(t\right), \control\left(t\right)\right)  \, dt 
\\
\\
\dot{\state}=\mathbf{f}\left(\state, \control\right), \quad \state(0) = \state_0,
\quad
\state\left( t_f \right) \in \target
\end{array}
\end{equation}
\noindent For optimal control problems with a continuous value function, the value function $V: \mathbb{R}^n \rightarrow [0, \infty)$ is the unique viscosity solution to the Hamilton-Jacobi-Bellman (HJB) equation. However, in many cases the value function may be non-differentiable or discontinuous at or near the target and the notion of envelope solution, as in \cite{bardi}, is needed to identify the value function as the solution to the HJB equation. Solving a partial differential equation is difficult in general because of the curse of dimensionality and here it is complicated further by looking for the solution in the viscosity sense.

In the literature, optimal feedback control for this problem has been done via discretization then applying a fixed point operation over a grid, mesh, or tensor decomposition of a grid/mesh as in \cite{bardi, cristiani_initialization_2010,  falcone_semi_lagrangian_2014, gorodetsky_2018}. The classic trade-off for mesh based methods is although accuracy is very high, scaling to dimensions higher than three is severely limited due to the requirement of adding exponentially more grid points per dimension. Some tensor decomposition methods alleviate the curse of dimensionality while retaining high accuracy so long as the tensor-rank of the optimal value function is low. That is, they are capable of trading dimension dependence with an auxillary metric of complexity, but the tensor-rank of arbitrary value functions may not always be low (even in low dimensions) and is hard to verify before performing the computation and observing the rank or time complexity. We aim to reduce dimension/tensor-rank dependence by applying a neural network based Reinforcement Learning approach. In order to do this, we show that it is possible to apply Reinforcement Learning to free terminal time type problems via the Semi-Lagrangian discretization, allowing one to convert the problem of interest (1) into a form resembling an infinite horizon optimal control problem with a varying discount.

\subsection{Semi-Lagrangian Discretization}
The Semi-Lagrangian method has been successfully used to solve low dimensional feedback control problems via mesh-based solvers for the HJB equation such as in the appendix of \cite{bardi, sl_book} or by seeding an open-loop method as in \cite{cristiani_initialization_2010}. We will apply the Kru\v{z}kov transformation, a classical transformation for minimum-time type problems,
\begin{equation}
\tilde{V}\left(\state\right) = \begin{cases}
    1 - e^{-V(\state)}, & \text{if $V(\state)<\infty$}.\\
    1, & \text{otherwise}.
  \end{cases}
\end{equation}
as is done in \cite{bardi, cristiani_initialization_2010} in order to acquire a functional equation which we will iteratively fit a solution to using a neural network.
The main idea will be to discretize the dynamics, in this case with forward Euler, to get a discrete optimal control problem which converges to the original problem as the time step $\Delta t$ decreases to zero. The value function approximation corresponding to the discrete dynamics should satisfy the discrete time dynamic programming principle. Therefore, instead of discretizing the directional derivative via finite differencing, we directly apply the dynamic programming principle for the discrete dynamics in combination with the same order integration method over $[0, \Delta t]$ as in \cite{falcone_semi_lagrangian_2014}:
\begin{align}
& \left \{ \begin{array}{l}
\state_{k+1} = \state_k + \Delta t \, \mathbf{f}(\state_k, \control)\\
\\
\int_{0}^{\Delta t} \ell(\state, \control) \, dt \approx \Delta t \, \ell(\state, \control)
\\
\end{array} \right .
\\
\implies &
V(\state_k) = \min_{\control \in \mathcal{U}}\left\{ \Delta t \, \ell(\state_k, \control) + V\left(\state_{k+1} \right) \right\}
\end{align}

Upon manipulation of both sides to rewrite the equation in terms of the Kru\v zkov transform of the value function we obtain
\[
\small
\begin{array}{ccc}
e^{-V(\state_k)} &=& e^{-\min_{\control} \left\{ \Delta t \ell(\state_k,\control) + V(\state_{k+1}) \right\} }\\[2mm]
e^{-V(\state_k)} &=& \underset{\control}{\max} \, e^{-\Delta t \ell(\state_k,\control) - V(\state_{k+1})}\\[2mm]
1-e^{-V(\state_k)} &=& 1-\underset{\control}{\max} \, e^{-\Delta t \ell(\state_k,\control) - V(\state_{k+1})}\\[2mm]
\tilde{V}(\state_k) &=& 1 + \underset{\control}{\min} \left\{ e^{-\Delta t \ell(\state_k,\control)}\left(\tilde{V}(\state_{k+1}) - 1\right)\right\}
\end{array}
\]

\noindent Furthermore, define the discount factor and operator for the right-hand-side of the (semi-discrete) HJB equation as
\begin{equation}
\begin{array}{cc}
\gamma(\state,\control) = e^{-\Delta t \ell(\state_k,\control)}\\
\tilde{H}\left( \state, \tilde{V}, \control \right) = 1 + \gamma\left(\state, \control\right)\left(\tilde{V}\left(\state + \Delta t \, \mathbf{f}\left(\state, \control\right) \right) - 1\right)
\end{array}
\end{equation}
\noindent as well as a bounded subset of $\mathbb{R}^n \setminus \target = \Omega$ to be the computational domain with which to solve over. Then the equation we wish to solve is

\begin{equation}\label{SL_HJB2}
\begin{array}{cc}

\tilde{V}(\state) = \underset{\control}{\min}\,\tilde{H}\left( \state, \tilde{V}, \control \right), & \state \in \Omega
\\
\tilde{V}(\state) = 0, & \state \in \target
\end{array}
\end{equation}

\noindent To obtain lower truncation error, one can use a second order differential equation discretization and integration scheme such as second order Runge-Kutta and the Trapezoid Rule. Within this paper, we choose to use Euler discretization in anticipation of high neural network fitting error dominating the overall error in computations and the complexity of backpropagation through the integration terms for $\state_{k+1}$ and the running cost. We formulate the problem in this way in order to force a discount factor, prevent gradients from appearing in the semi-discrete HJB equation we wish to solve, and allow for a surrogate model for both value and control functions. Most importantly, if the discount factor $\gamma$ was a fixed number then this equation would correspond to an infinite horizon discrete optimal control problem where $\gamma$ is used to make the running cost finite when integrated over an infinite horizon. Such a problem is amenable to approximate dynamic programming methods as studied in \cite{bertsekas_neurodynamic} as well as reinforcement learning methods as in \cite{lewis_reinforcement_2009, tutsoy_reinforcement_2016}. In some cases, such as in \cite{tutsoy_reinforcement_2016}, the minimum time problem is approximated by setting $\gamma \approx 1$ and applying reinforcement learning techniques such as value iteration or an Actor-Critic method. However, rate of convergence and the conditioning of the problem become worse as $\gamma \rightarrow 1$ as the step size and optimal cost are directly related to $\gamma$. We aim to make use of this interpretation in order to apply the well known Actor-Critic paradigm, see \cite{lewis_reinforcement_2009, scherrer_approximate_2015, lillicrap, zhou_actor-critic_2021}, which is typically applied to infinite horizon optimal control problems.

\section{Actor-Critic Algorithm for Kru\v{z}kov HJB}
Motivated by the success of Actor-Critic methods for infinite horizon optimal control, we propose to use a Least Squares Temporal Difference Actor-Critic framework in the case studied here. We implement multilayer neural networks to represent both the actor (controller) $U(\state, W_U)$ and critic (value function) representation, $V(\state, W_V)$. We do this to mitigate the necessity of both choosing good set of basis functions for approximating the true optimal control/cost as well as having enough of such functions to represent the value and control functions adequately. It should be noted that although this is effective in practice, the compatibility of the actor and critic will no longer be guaranteed and the gradients $\partial_\control \tilde{H}(\state, \tilde{V}, \control)$ may be biased estimates of the true descent direction (\cite{silver_deterministic_2014}). 

The actor critic algorithm presented here in Algorithm \ref{algo:actor_critic} is based on value iteration, where the control and value approximators are optimized in an alternating fashion. Upon each iteration $k$, the critic $\tilde{V}(\state, W_V)$ is initially held fixed while the control approximation $U(\state, W_U)$ minimizes its weights against the expected cost-to-go of one step in time, $\tilde{H}(\state, \tilde{V}, \control)$ obtaining $W_U^{k+1}$. After the control refinement stage, the new value approximation is set to the closest approximation to the new predicted cost-to-go, $\tilde{H}(\state, \tilde{V}, U(\state, W_U^{k+1})$. If this approximation is near enough, we expect both $U(\state, W_U)$ and $V(\state, W_V)$ to converge to the optimal control and value function respectively. This intuition is due to the fact that if we use meshes coupled with the Markov Chain Approximation of \cite{kushner_dupuis} then provided $1 \leq \ell(\state, \control) \leq M$ on $\Omega$ and some technical assumptions on the dynamics and target set as outlined in \cite{bardi} then this process would converge to the desired approximations.

\begingroup
\setlength{\algomargin}{5mm}
\begin{algorithm2e}
\DontPrintSemicolon 
\KwIn{A tuple of weights $(W_V^0, W_U^0)$ for the value and actor networks, a time step $\Delta t$, a moving average coefficient $\alpha \in (0, 1]$, and the number of samples to draw from each set $N_\Omega, N_\target$.}
\KwOut{Neural network approximations of cost-to-go and optimal control.}
\vspace{2mm}
\small
Let $\widehat{W}_V \leftarrow W_V^0$\;
Let $\widehat{W}_U \leftarrow W_U^0$\;

\vspace{2mm}
\While{$\widehat{W}_V, \widehat{W}_U$ not converged}{
\vspace{2mm}
Draw set of samples $X_\Omega \leftarrow \left\{ \state_i \right\}_{i=1}^{N_\Omega} \sim \mathbb{P}(\Omega$)\;
Draw set of samples $X_\target \leftarrow \left\{ \state_j \right\}_{j=1}^{N_\target} \sim \mathbb{P}(\target$)\;

Define function $\tilde{H}\left(\state, W \right) := \tilde{H}\left(\state, \tilde{V}(\state, \widehat{W}_V), U(\state, W)\right)$\;

\vspace{2mm}
$W_U^{k+1} \leftarrow \underset{W_U}{\text{argmin}} \underset{\state \in X_\Omega}{\mathbb{E}} \left[ \tilde{H}\left(\state, W_U \right) \right]$\;

\vspace{2mm}
$\widehat{W}_U \leftarrow \alpha W^{k+1}_U + (1 - \alpha) \widehat{W}_U$\;

\vspace{4mm}$L^k_{X_\Omega}(W_V) \leftarrow \underset{\state \in X_\Omega}{\mathbb{E}}\left[\left(\tilde{V}(\state, \widehat{W}_V) - \tilde{H}(\state, \widehat{W}_U)\right)^2\right]$\;
\vspace{1mm}
$L^k_{X_\target}(W_V) \leftarrow\underset{\state \in X_\target}{\mathbb{E}} \left[\tilde{V}\left(\state, W_V\right)^2\right] $\;

\vspace{4mm}
$W_V^{k+1} \leftarrow \underset{W_V}{\text{argmin}} \, \left[ L^k_{X_\Omega}(W_V) + L^k_{X_\target}(W_V) \right]$\;

\vspace{4mm}
$\widehat{W}_V \leftarrow \alpha W^{k+1}_V + (1 - \alpha) \widehat{W}_V$\;
}
\Return{ $(\widehat{W}_V, \widehat{W}_U)$}
\caption{{\sc Actor-Critic Iteration}}
\label{algo:actor_critic}
\end{algorithm2e}
\endgroup

\noindent This algorithm can be efficiently performed in parallel over the states because neural networks can operate on batch inputs efficiently and mature frameworks for their batched gradient computation are readily available. Practically, we implement this algorithm as a stochastic gradient method, only taking a few steps per optimization stage (stages 7 and 11) using an optimizer such as Adam from \cite{adam_2017} and increasing the number of steps at stage 11 if the fitting error is large. Then we update a moving average of the model weights $\widehat{W}_V, \widehat{W}_U$to reduce variance in the weight updates. Algorithm (1) can be interpreted as an Actor-Critic type fitted value iteration (or policy gradient when steps 7 and 11 are single gradient descent updates) where each sample from the domain is propagated forward in time by $\Delta t$ through an actor $U(\state, W_U)$ and the new state is used to update the critic's value estimate. If too few samples are drawn, weight updates will have a high amount of noise and the algorithm may become unstable. One can reduce the update noise through the use of a target network as in \cite{lillicrap} where the weights are slowly averaged into an auxillary network in order to provide stable learning dynamics. We find that this is not necessary for our two examples where the value function is continuous, but helps reduce weight oscillation when there is a discontinuity. Algorithm \ref{algo:actor_critic} may be able to be transformed into a Monte Carlo, rollout-based scheme as in \cite{zhou_actor-critic_2021}, but due to deterministic dynamics we are able to take advantage of per-sample parallelism without the need for propagation beyond one step. We also note that in the stochastic case, conversion to the form of equation \ref{SL_HJB2} should not be necessary since the system's diffusion should guarantee a unique, sufficiently differentiable solution to the original HJB equation.

\begin{remark}
This algorithm can be extended to use supervised data as in \cite{nakamura-zimmerer_adaptive_2021} by generating a set of known optimal tuples  and appending a loss to steps 7 or 11.
\end{remark}

\begin{remark}
A more efficient domain sampling which takes into account model error over the domain instead of uniform random sampling would likely lead to faster convergence. In regards to convergence criteria and iteration count, either a fixed number of iterations is done and the models are evaluated by closed-loop performance or the mean and variance of the moving averages $\widehat{W}_V, \widehat{W}_U$ are used to stop the iteration.
\end{remark}

\section{Examples}

To demonstrate the algorithm described, we apply it to three qualitatively different examples. The first is the double integrator, a canonical academic example known to many in optimal control, which demonstrates a discontinuous (bang-bang) optimal control and non-differentiable value function.. The second is the control-regularized Dubins Vehicle which displays a value function which is both non-differentiable and discontinuous in some regions of its state-space. These properties make optimal open-loop trajectories difficult to obtain without expert initial guesses as there are infinitely many locally optimal controls. The last example is a minimum time attitude control problem for a 7-dimensional rigid body to demonstrate scalability of the algorithm. For each example, we apply the Kru\v{z}kov transformation to the critic's final output during training as well as enforce control constraints via the actor's final layer.

\subsection{Double Integrator}
Our first example is that of the time-optimal double integrator:
\begin{equation}
\begin{array}{c}
\underset{\control, t_f}{\min} \int_{0}^{t_f} dt
\\
\begin{bmatrix}
\dot{x} \\
\dot{y} 
\end{bmatrix}
=
\begin{bmatrix}
y \\
u
\end{bmatrix}, \quad
x(t_f) = y(t_f) = 0, \quad \vert u \vert \leq 1
\end{array}
\end{equation}

\noindent This is a simple example of a problem which gives rise to a non-smooth value function without discontinuity. The switching curve is known to be parabolic and the control is piecewise $+1$ or $-1$ on each side of this curve, given by $\control(\state) = -\text{sign}(D_y V(\state))$ wherever $V$ is differentiable. However, generating optimal trajectories for this system via open-loop methods presents difficulties due to the discontinuous optimal control. For computation, we aim to solve this problem over the ball $B(r=5)$ and apply the actor-critic algorithm with networks of 4 layers each comprised of 128 neurons, a time discretization of $dt = 0.05$, no moving average ($\alpha=0$), and sample sizes of $\vert X_\Omega \vert = 500, \, \vert X_{\target} \vert = 1$.

\begin{figure}[!htb]
\centering
\includegraphics[scale=0.37]{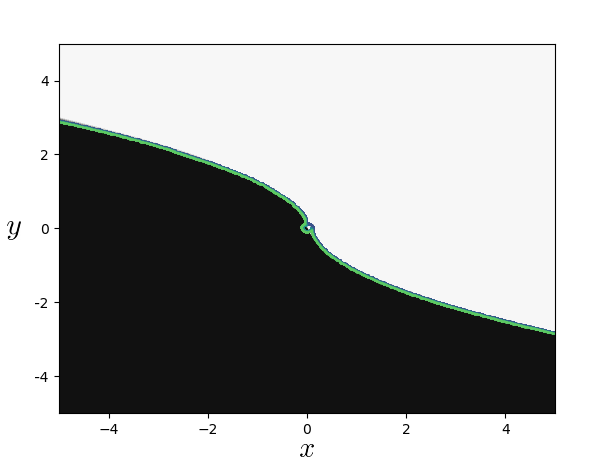}
\caption{Switching surface given by neural network actor with tanh activation (white/black) and switching curve computed by grid-based semi-lagrangian iteration (green)}
\label{fig:integrator-control}
\end{figure}

We use tanh activations in this case because the optimal control is bang-bang and the tanh function can smoothly saturate between the bounds +1 and -1. Notably, the actor network obtained is continuous by construction and therefore the switching behavior is approximated by a steep, smooth transition between -1 to +1 (and vice-versa) which also happens to avoid the chattering problem faced by bang-bang controllers.

\FloatBarrier

\begin{figure}[!hb]
  \centering
\includegraphics[scale=0.37]{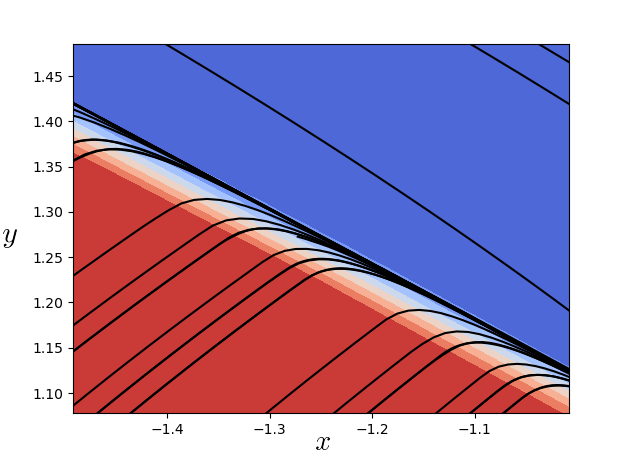}
\caption{Close up of switching curve given by neural network actor with sample trajectories overlayed (black). }
\label{fig:integrator-switching}
\end{figure}

\begin{figure}[!ht]
\centering
\includegraphics[scale=0.75]{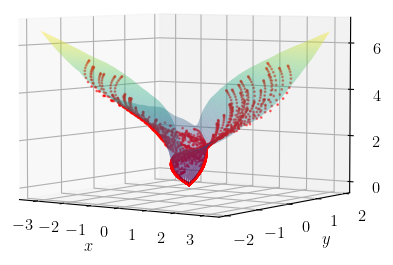}
\caption{Comparison between the true value function (red) and the critic network (surface).}
\label{fig:integrator-3dplot-comparison}
\end{figure}
\FloatBarrier
Compared to the true solution over a grid of 3200 test points, obtained using a pseudospectral collocation method, see \cite{ross_review_2012}. Note the critic (value approximator) network is very similar geometrically to the optimal value function and has typically differs from the target with mean-squared-error of order $10^{-2}$.
\subsection{Dubins Vehicle}
For our final example we will consider the control-regularized minimum time problem for Dubin's Vehicle:
\begin{equation}
\begin{array}{c}
\underset{u, t_f}{\min} \int_{0}^{t_f} \left( 1+\frac{1}{2}\left[u(t)\right]^2 \right)  \, dt 
\\
\\
\begin{bmatrix}
\dot{x} \\
\dot{y} \\
\dot{\theta}
\end{bmatrix}
=
\begin{bmatrix}
\cos{\theta}\\
\sin{\theta}\\
u
\end{bmatrix}
, \quad
\begin{array}{c}
\vert u \vert \leq 6, \quad t_f \geq 0\\
\vert x(t_f) \vert^2 + \vert y(t_f) \vert^2 \leq \left( 0.1 \right)^2
\end{array}
\end{array}
\end{equation}

\noindent This problem displays a value function with both a non-smooth gradient and a discontinuity at the boundary of the target set. Physically, we aim to steer a vehicle with constant velocity to the origin in the least time and using the least steering input. The discontinuity at the boundary of the target set comes from the fact that the dynamics are not small-time-locally-controllable \cite{bardi}. In other words, starting near the target does not guarantee a small cost-to-go for some states because a lack the controllability leads to a high cost-to-go as the vehicle must follow an arc dependent on its turning radius. For example, starting at position $\left(x, y, \theta \right) = \left(0.1, 0, 0\right)$ the vehicle must first travel away from $\target$ and circle back since there is no direct control of $x, y$ and the velocity of the vehicle is constant. This problem has a radially symmetric value function and so visualization is done in only the $(x, 0, \theta)$ coordinates. For reference, we present below a plot of a representative slice of the value function as computed via the Markov Chain Approximation of Kushner and Dupuis \cite{kushner_dupuis}.

\begin{figure}[!htb]
  \centering
  \includegraphics[scale=0.5]{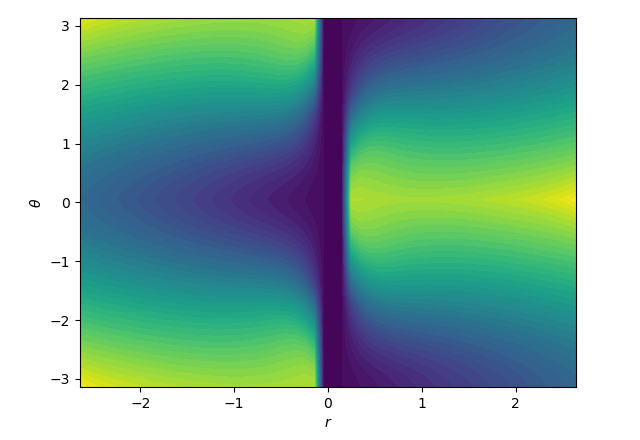}
  \caption{The reference value function as computed via mesh as in \cite{kushner_dupuis}, shown radially in the $(r, \theta)$ plane where $r$ is on the x-axis and $\theta$ on the y-axis.}
  \label{fig:reference}
\end{figure}
\begin{figure}[!htb]
\centering
\includegraphics[scale=0.5]{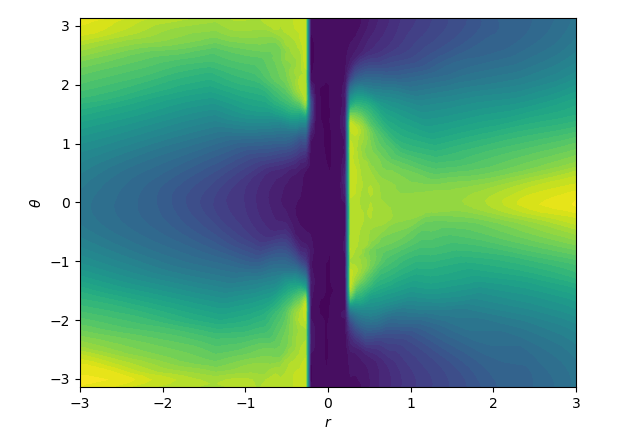}
\caption{Contour of the neural network approximation of the value function in the  $(r, \theta)$ plane.}\label{fig:dubins-value}
\end{figure}

The Actor-Critic method using a weight-averaging update for the target network with $\alpha = 0.1$, ReLU activation with 4 layers of 128 neurons per layer, and a rescaling of the value function by $0.2$ produced excellent results. Although the optimal control is discontinuous, owing to the fact that the vehicle must physically turn left or right depending on heading, the neural network controller $U(\state, W_U)$ performs well after tuning the parameter $\alpha$ governing the moving average.

\begin{figure}[!htb]
  \centering
  \includegraphics[scale=0.35]{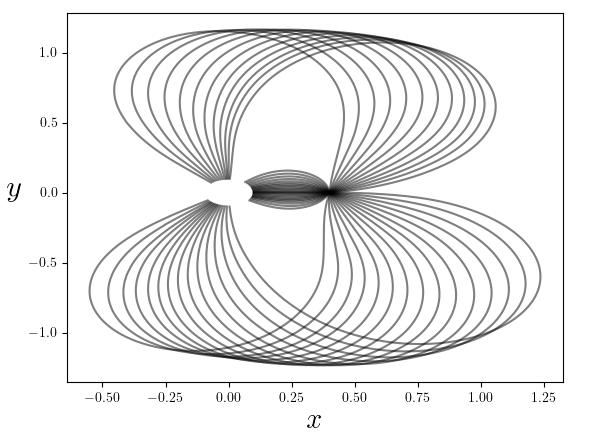}
  \caption{Actor feedback control over sample trajectories in the $(x, y)$ space near the target. Each trajectory begins with an initial heading $\theta(0) \in [-\pi, \pi)$. There is a paradigm switch between trajectories facing near the target to begin to face about 45 degrees away, corresponding to a kink in the value function. The choice of steering left or right when facing exactly opposite to the target also corresponds to a kink.}
  \label{fig:trajectories}
\end{figure}

We found that it is necessary to rescale either the running cost $\ell(\state, \control)$ in order to produce adequate approximations for this problem as the exponential term $\exp(-V(\state)) \approx 0$ when $V(\state) \approx 5$ and many states do achieve a cost greater than that limit. Due to limited precision of floating point, the Kru\v{z}kov transformation truncates the values of all states which would incur a cost of about 5. We currently get around this limitation through either the rescaling $\ell(\state, \control) \rightarrow \mu \, \ell(\state, \control)$ or a direct rescaling of the dynamic programming equation through $V(\state) \rightarrow \mu \, V(\state)$ for $0 < \mu \leq 1$.

\subsection{Attitude Control of Spacecraft}
For our next example, we chose a higher dimensional satellite system to show that this algorithm has the potential to push to higher dimensions. The optimal control problem is given as
\begin{equation}
\begin{array}{c}

\underset{u, t_f}{\min} \quad
t_f
\\
\begin{bmatrix}
\dot{q}_0\\
\dot{\mathbf{q}}\\
\mathbf{J}\dot{\boldsymbol{\omega}}
\end{bmatrix}
=
\begin{bmatrix}
-\frac{1}{2}\boldsymbol{\omega}^T \mathbf{q}\\
\frac{1}{2}\left( \mathbf{-\omega}_\times \mathbf{q} + q_0 \boldsymbol{\omega} \right)\\
-\boldsymbol{\omega}_\times \mathbf{J} \boldsymbol{\omega} - \control
\end{bmatrix}
\\
\\
-0.3 \leq \mathbf{u} \leq 0.3\\
\mathbf{q}(t_f) = \boldsymbol{\omega}(t_f) = 0, \quad q_0(t_f) = 1
\\
\\
J = \begin{bmatrix}
59.22 & -1.14 & -0.8 \\
-1.14 & 40.56 & 0.1 \\
-0.8 & 0.1 & 57.60
\end{bmatrix}
\end{array}
\end{equation}
The state space sampling is done by sampling Euler angles in $[-\frac{\pi}{2}, \frac{\pi}{2}]$ uniformly then converting to quaternions, the angular velocities are sampled uniformly from the sphere $10^{-4} \leq \Vert \boldsymbol{\omega}\Vert^2 \leq 0.3$, and the inertia matrix is that of the TRACE spacecraft (see \cite{trace}). We used the Adagrad optimizer, sample sizes of $\vert X_\Omega \vert = \vert X_\target \vert = 1500$, $\Delta t = 0.3$, and ran the algorithm for 2000 iterations for a total of 73 seconds. The large $\Delta t$ is used to speed up convergence as the time-scale of trajectories is comparatively large. The neural networks used are sequential networks with hidden layers, $L_i(x)$, of sizes $[7 \rightarrow 200 \rightarrow 200]$ where the hidden layers are residual blocks of the form $x \rightarrow L_{i-1}(x) + x$, and ReLU activation.

\begin{figure}[htp]
  \centering
  \includegraphics[scale=0.5]{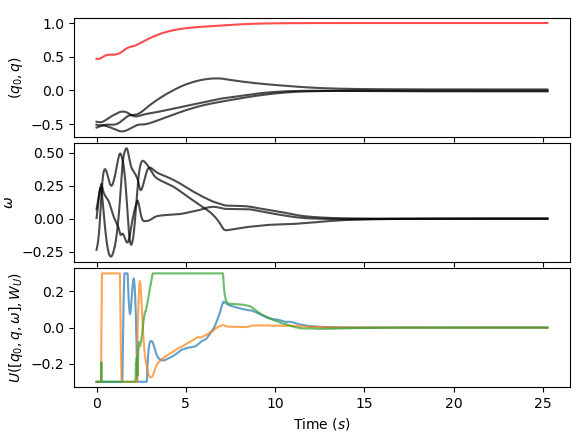}
  \caption{An example trajectory for the TRACE system using the actor network as the controller. The value of $q_0$ is in red and $\mathbf{q}, \boldsymbol{\omega}$ in black. Note that the control inputs resemble that of a smoothed, bang-bang controller which one expects to observe for minimum time solutions.}
  \label{fig:reference3}
\end{figure}

\begin{figure}[htp]
  \centering
  \includegraphics[scale=0.5]{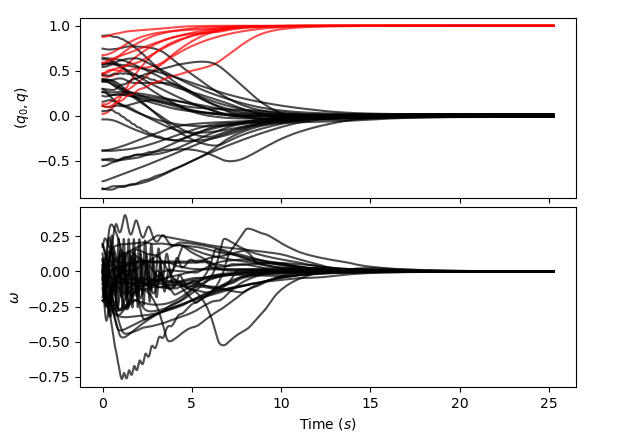}
  \caption{Ensemble trajectories displaying stability of the neural network controller.}
  \label{fig:reference2}
\end{figure}

It is important to note here that because of the generality of the model, the terminal condition $\Vert \mathbf{q} \Vert = 0$ may not exactly correspond to $\tilde{V}(q(t_f), W_V) = 0$ and therefore the steady-state may be $\Vert \mathbf{q} \Vert = \epsilon$ since the model is minimizing its squared deviation from 0 over the target. Experimentally, we have observed that for nearly all cases that the steady state satisfies $\Vert \mathbf{q} \Vert_\infty \leq 0.05$. Modifications to achieve neural network controller stability can be found \cite{nakamura-zimmerer_neural_2022}. Updating the model architecture to automatically accommodate the terminal condition would be ideal and enforce $\mathbf{q}(t_f) = 0$ exactly. A collection of 100 sample trajectories are shown in figure \ref{fig:reference2}, note that they are stable. In future work, we intend to compare these trajectories with open-loop solutions generated via pseudospectral methods.

\section{Conclusion}
In this paper we demonstrated the use of a neural network Actor-Critic algorithm for solving optimal control problems with free terminal time that may have nonsmooth solutions in the value or control. We find that it is a promising approach towards alleviating the curse of dimensionality and as a possible companion method to open-loop solvers that struggle with local optima or computing a solution near kinks and corners of the value function. We find that the method has promise even for problems with discontinuous optimal cost, such as in the Dubins Vehicle example. In future work, the efficacy of coupling the proposed algorithm to direct methods could be interesting as well as analysis on the convergence of such methods.

\section*{Acknowledgement}
This work was supported by the Air Force Office of Scientific Research (AFOSR) under grant no. FA9550-21-1-0113, and the National Science Foundation (NSF) under grant no. 2134235.

\bibliographystyle{unsrtnat}
\bibliography{ifacconf}

\begin{thebibliography}{18}
\providecommand{\natexlab}[1]{#1}
\providecommand{\url}[1]{\texttt{#1}}
\expandafter\ifx\csname urlstyle\endcsname\relax
  \providecommand{\doi}[1]{doi: #1}\else
  \providecommand{\doi}{doi: \begingroup \urlstyle{rm}\Url}\fi

\bibitem[Bardi and Dolcetta(1997)]{bardi}
Martino Bardi and Italo~Capuzzo Dolcetta.
\newblock \emph{Optimal {Control} and {Viscosity} {Solutions} of
  {Hamilton}-{Jacobi}-{Bellman} {Equations}}.
\newblock Modern {Birkhauser} {Classics}. Birkhauser Basel, 1997.
\newblock ISBN 978-0-8176-4754-4.
\newblock \doi{10.1007/978-0-8176-4755-1}.

\bibitem[Cristiani and Martinon(2010)]{cristiani_initialization_2010}
E.~Cristiani and P.~Martinon.
\newblock Initialization of the {Shooting} {Method} via the
  {Hamilton}-{Jacobi}-{Bellman} {Approach}.
\newblock \emph{Journal of Optimization Theory and Applications}, 146\penalty0
  (2):\penalty0 321--346, August 2010.
\newblock ISSN 1573-2878.
\newblock \doi{10.1007/s10957-010-9649-6}.

\bibitem[Falcone et~al.(2014)Falcone, Kalise, and
  Kröner]{falcone_semi_lagrangian_2014}
Maurizio Falcone, Dante Kalise, and Axel Kröner.
\newblock A semi-{Lagrangian} scheme for {Lp}-penalized minimum time problems.
\newblock In \emph{21st {International} {Symposium} on {Mathematical} {Theory}
  of {Networks} and {Systems}}, Groningen, Netherlands, July 2014.

\bibitem[Gorodetsky et~al.(2018)Gorodetsky, Karaman, and
  Marzouk]{gorodetsky_2018}
Alex Gorodetsky, Sertac Karaman, and Youssef Marzouk.
\newblock High-dimensional stochastic optimal control using continuous tensor
  decompositions.
\newblock \emph{The International Journal of Robotics Research}, 37\penalty0
  (2-3):\penalty0 340--377, February 2018.
\newblock ISSN 0278-3649.
\newblock \doi{10.1177/0278364917753994}.
\newblock Publisher: SAGE Publications Ltd STM.

\bibitem[Falcone and Ferretti(2013)]{sl_book}
Maurizio Falcone and Roberto Ferretti.
\newblock \emph{Semi-{Lagrangian} {Approximation} {Schemes} for {Linear} and
  {Hamilton} {Jacobi} {Equations}}.
\newblock Other {Titles} in {Applied} {Mathematics}. Society for Industrial and
  Applied Mathematics, December 2013.
\newblock ISBN 978-1-61197-304-4.
\newblock \doi{10.1137/1.9781611973051}.

\bibitem[Bertsekas and Tsitsiklis(1996)]{bertsekas_neurodynamic}
Dimitri~P. Bertsekas and John~N. Tsitsiklis.
\newblock \emph{Neuro-{Dynamic} {Programming}}.
\newblock Athena Scientific, 1st edition, 1996.
\newblock ISBN 978-1-886529-10-6.

\bibitem[Lewis and Vrabie(2009)]{lewis_reinforcement_2009}
Frank~L. Lewis and Draguna Vrabie.
\newblock Reinforcement learning and adaptive dynamic programming for feedback
  control.
\newblock \emph{IEEE Circuits and Systems Magazine}, 9\penalty0 (3):\penalty0
  32--50, 2009.
\newblock ISSN 1558-0830.
\newblock \doi{10.1109/MCAS.2009.933854}.
\newblock Conference Name: IEEE Circuits and Systems Magazine.

\bibitem[Tutsoy and Brown(2016)]{tutsoy_reinforcement_2016}
Onder Tutsoy and Martin Brown.
\newblock Reinforcement learning analysis for a minimum time balance problem.
\newblock \emph{Transactions of the Institute of Measurement and Control},
  38\penalty0 (10):\penalty0 1186--1200, October 2016.
\newblock ISSN 0142-3312.
\newblock \doi{10.1177/0142331215581638}.
\newblock Publisher: SAGE Publications Ltd STM.

\bibitem[Scherrer et~al.(2015)Scherrer, Ghavamzadeh, Gabillon, Lesner, and
  Geist]{scherrer_approximate_2015}
Bruno Scherrer, Mohammad Ghavamzadeh, Victor Gabillon, Boris Lesner, and
  Matthieu Geist.
\newblock Approximate {Modified} {Policy} {Iteration} and its {Application} to
  the {Game} of {Tetris}.
\newblock \emph{Journal of Machine Learning Research}, 16\penalty0
  (49):\penalty0 1629--1676, 2015.
\newblock ISSN 1533-7928.

\bibitem[Lillicrap et~al.(2016)Lillicrap, Hunt, Pritzel, Heess, Erez, Tassa,
  Silver, and Wierstra]{lillicrap}
Timothy~P. Lillicrap, Jonathan~J. Hunt, Alexander Pritzel, Nicolas Heess, Tom
  Erez, Yuval Tassa, David Silver, and Daan Wierstra.
\newblock Continuous control with deep reinforcement learning.
\newblock In Yoshua Bengio and Yann LeCun, editors, \emph{4th International
  Conference on Learning Representations, {ICLR} 2016, San Juan, Puerto Rico,
  May 2-4, 2016, Conference Track Proceedings}, 2016.

\bibitem[Zhou et~al.(2021)Zhou, Han, and Lu]{zhou_actor-critic_2021}
Mo~Zhou, Jiequn Han, and Jianfeng Lu.
\newblock Actor-{Critic} {Method} for {High} {Dimensional} {Static}
  {Hamilton}--{Jacobi}--{Bellman} {Partial} {Differential} {Equations} based on
  {Neural} {Networks}.
\newblock \emph{SIAM Journal on Scientific Computing}, 43\penalty0
  (6):\penalty0 A4043--A4066, January 2021.
\newblock ISSN 1064-8275.
\newblock \doi{10.1137/21M1402303}.
\newblock Publisher: Society for Industrial and Applied Mathematics.

\bibitem[Silver et~al.(2014)Silver, Lever, Heess, Degris, Wierstra, and
  Riedmiller]{silver_deterministic_2014}
David Silver, Guy Lever, Nicolas Heess, Thomas Degris, Daan Wierstra, and
  Martin Riedmiller.
\newblock Deterministic {Policy} {Gradient} {Algorithms}.
\newblock In \emph{Proceedings of the 31st {International} {Conference} on
  {Machine} {Learning}}, pages 387--395. PMLR, January 2014.
\newblock ISSN: 1938-7228.

\bibitem[Kushner and Dupuis(2001)]{kushner_dupuis}
Harold Kushner and Paul~G. Dupuis.
\newblock \emph{Numerical {Methods} for {Stochastic} {Control} {Problems} in
  {Continuous} {Time}}.
\newblock Stochastic {Modelling} and {Applied} {Probability}. Springer-Verlag,
  New York, 2 edition, 2001.
\newblock ISBN 978-0-387-95139-3.
\newblock \doi{10.1007/978-1-4613-0007-6}.

\bibitem[Kingma and Ba(2017)]{adam_2017}
Diederik~P. Kingma and Jimmy Ba.
\newblock Adam: A method for stochastic optimization, 2017.

\bibitem[Nakamura-Zimmerer et~al.(2021)Nakamura-Zimmerer, Gong, and
  Kang]{nakamura-zimmerer_adaptive_2021}
Tenavi Nakamura-Zimmerer, Qi~Gong, and Wei Kang.
\newblock Adaptive {Deep} {Learning} for {High}-{Dimensional}
  {Hamilton}–{Jacobi}–{Bellman} {Equations}.
\newblock \emph{SIAM Journal on Scientific Computing}, 43\penalty0
  (2):\penalty0 A1221--A1247, 2021.
\newblock Publisher: SIAM.

\bibitem[Ross and Karpenko(2012)]{ross_review_2012}
I.~Michael Ross and Mark Karpenko.
\newblock A review of pseudospectral optimal control: {From} theory to flight.
\newblock \emph{Annual Reviews in Control}, 36\penalty0 (2):\penalty0 182--197,
  December 2012.
\newblock ISSN 1367-5788.
\newblock \doi{10.1016/j.arcontrol.2012.09.002}.

\bibitem[Zimbelman et~al.(1995)Zimbelman, Wilmot, and Evangelista]{trace}
Darrell Zimbelman, Jonathan Wilmot, and Solomon Evangelista.
\newblock The {Attitude} {Control} {System} {Design} for the {Transition}
  {Region} and {Coronal} {Explorer} {Mission}.
\newblock \emph{Small Satellite Conference}, September 1995.

\bibitem[Nakamura-Zimmerer et~al.(2022)Nakamura-Zimmerer, Gong, and
  Kang]{nakamura-zimmerer_neural_2022}
Tenavi Nakamura-Zimmerer, Qi~Gong, and Wei Kang.
\newblock Neural {Network} {Optimal} {Feedback} {Control} with {Guaranteed}
  {Local} {Stability}, May 2022.
\newblock arXiv:2205.00394 [cs, eess, math].

\end{thebibliography}

\end{document}